\documentclass[12 pt]{amsart}

\usepackage{graphicx}
\usepackage{color}
\usepackage{amssymb, amsmath, amsthm}
\usepackage{mathptmx}

\usepackage{epsfig}
\usepackage{epstopdf}
\usepackage{graphicx}

\theoremstyle{plain}

\newtheorem{theorem}{Theorem}
\newtheorem{lemma}{Lemma}
\newtheorem{corollary}{Corollary}

\newtheorem*{mtheorem}{Main Theorem}
\newtheorem{conjecture}{Conjecture}

\newcommand{\bi}{\begin{itemize}}
\newcommand{\ei}{\end{itemize}}
\newcommand{\be}{\begin{enumerate}}
\newcommand{\ee}{\end{enumerate}}

\newcommand{\n}{\beta}

\newcommand{\emp}{\emptyset}
\newcommand{\X}{\times}

\newcommand{\eps}{\epsilon}

\newcommand{\A}{\alpha}
\newcommand{\pd}{\partial}

\numberwithin{definition}{section}
\numberwithin{example}{section}
\numberwithin{lemma}{section}
\numberwithin{theorem}{section}
\numberwithin{corollary}{section}

\begin{document}
\title{Uniqueness of higher genus bridge surfaces for torus knots}
\author{Alexander Zupan}
\thanks{The author is supported by the National Science Foundation under Award No. DMS-1203988.}

\maketitle

\begin{abstract}
We show that a torus knot which is not 2-bridge has a unique irreducible bridge splitting of positive genus.
\end{abstract}

\section{Introduction}

A bridge splitting of a knot $K$ in a compact orientable 3-manifold $M$ decomposes the pair $(M,K)$ into two simple pieces attached along a bridge surface.  As such, a bridge splitting may be viewed as a direct extension of the notion of a Heegaard splitting of a 3-manifold.  In keeping with the current direction of the theory of Heegaard splittings, mounting evidence connects the complexity of a bridge surface to the topology and geometry of the knot exterior.  Various measures of bridge surface complexity have been used to characterize the essential surfaces in a knot exterior \cite{bs}, to provide obstructions to a knot admitting some exceptional Dehn surgeries \cite{bbl, bcjtt}, and to estimate the hyperbolic volume of certain knots \cite{zupan1}. \\

A natural problem in this setting is to classify all of the bridge splittings for large collections of knots, up to the generic operations of perturbation, stabilization, and meridional stabilization.  A bridge surface is said to be \emph{irreducible} if it is not the result of one of these three operations.  One might expect that the simplest knots have the simplest bridge structure; namely, a knot whose topology is uncomplicated should have very few irreducible bridge surfaces.  Historically, Otal examined this idea for bridge spheres, proving that the unknot and each 2-bridge knot have a unique irreducible bridge sphere \cite{otal}.  More recently, Hayashi and Shimokawa have extended Otal's result to show that the 1-bridge sphere for the unknot is, not surprisingly, the unique irreducible bridge surface for the unknot \cite{hayshi}, and Scharlemann and Tomova prove a similar statement for 2-bridge knots \cite{schartom}.  In contrast, Ozawa and Takao have produced a knot with two irreducible bridge spheres of different bridge number \cite{ozawatakao}, and Jang has exhibited 3-bridge links with infinitely many distinct minimal bridge spheres \cite{jang}. \\

Torus knots provide another class of knots worthy of study in this regard.  Ozawa has shown that, like the unknot and 2-bridge knots, torus knots have a unique irreducible bridge sphere \cite{ozawa}.  In this paper, we extend his result to prove the following: \\

\begin{mtheorem}
Let $K$ be a torus knot which is not 2-bridge.  Then $K$ has a unique irreducible bridge surface of positive genus, a $(1,1)$-surface. \\
\end{mtheorem}

The proof of the theorem uses the fact that irreducible bridge surfaces may be partitioned into classes of strongly irreducible surfaces, weakly reducible surfaces, and a third class we call \emph{essentially cancelable}.  In Section \ref{strongly}, we show that a torus knot has a unique strongly irreducible bridge surface; namely, a $(0,\min\{p,q\})$-surface.  The argument relies on machinery developed by the author in \cite{zupan2} to study the decomposition of a strongly irreducible bridge surface induced by cutting along an essential surface.  In Section \ref{weak}, we show that every weakly reducible bridge surface for a torus knot is reducible by invoking results of Tomova regarding weakly reducible bridge surfaces \cite{tomova} and work of Boileau, Rost, and Zieschang \cite{brz}, Moriah \cite{moriah}, and Moriah and Schultens \cite{morsch} classifying Heegaard splittings of torus knot exteriors.  In Section \ref{esscan}, we examine the third class of surfaces and demonstrate that a non-2-bridge torus knot has a unique essentially cancelable surface, the $(1,1)$-surface mentioned in the statement of the main theorem.  Finally, Section \ref{conje} includes some closing remarks and a conjecture generalizing the main theorem. \\

\noindent \textbf{Acknowledgements} \, I would like to thank Cameron Gordon, Yeonhee Jang, Tye Lidman, and Maggy Tomova for helpful conversations and insights over the course of this project.

\section{Preliminaries}\label{pre}

Throughout this paper, we will let $\eta( \cdot)$ and $N( \cdot)$ denote open and closed regular neighborhoods, respectively, in an ambient manifold that should be clear from context.  Given a knot $K$ in a closed, orientable, irreducible 3-manifold $M$, we let $E(K) = M \setminus \eta(K)$ denote the \emph{exterior} of $K$. \\

Consider a handlebody $V$ containing a collection of properly embedded arcs $\A$ such that $\A$ is isotopic into $\Sigma = \pd V$.  Such a collection of arcs is called \emph{trivial}.  By virtue of this isotopy, we may suppose that each arc $\A_i \subset \A$ cobounds an embedded disk $\Delta_i$ with an arc $\gamma_i \subset \Sigma$, where $\Delta_i \cap \A = \A_i$.  In addition, the collection $\{\Delta_i\}$ may be chosen to be pairwise disjoint, and we call each $\Delta_i$ a \emph{bridge disk}.  A \emph{compressing disk} for the pair $(V,\A)$ is a properly embedded disk $D$ such that $\pd D$ is an essential curve in $\Sigma \setminus \eta(\A)$ and $D \cap \A = \emp$.  Similarly, a \emph{cut disk} is a properly embedded disk $C$ such that $\pd C$ is an essential curve in $\Sigma \setminus \eta(\A)$, with $C$ intersecting $\A$ in a single point. \\

Let $K$ be a knot in $M$.  A $(g,b)$\emph{-bridge splitting} of $(M,K)$ is a decomposition
\[ (M,K) = (V,\A) \cup_{\Sigma} (W,\n),\]
where $(V,\A)$ and $(W,\n)$ are genus $g$ handlebodies each containing a collection of $b$ trivial arcs, $\Sigma = \pd V = \pd W$, and $K = \A \cup \n$.  In this context, a $b$-bridge sphere is a $(0,b)$-surface.  Note that $\Sigma$ is also a Heegaard surface for $M$, and to avoid confusion, we let $\Sigma_K$ represent the compact genus $g$ surface with $2b$ boundary components $\Sigma \setminus \eta(K)$.  Two bridge surfaces $\Sigma$ and $\Sigma'$ are considered equivalent if $\Sigma_K$ is isotopic to $\Sigma'_K$ in $E(K)$. \\

Much of the terminology and machinery used to characterize Heegaard splittings of 3-manifolds can be adapted and expanded to understand bridge surfaces for knots in 3-manifolds.  Suppose that $(M,K) = (V,\A) \cup_{\Sigma} (W,\n)$ is a $(g,b)$-bridge splitting.  If there are compressing disks $D$ and $D'$ for $(V,\A)$ and $(W,\n)$ such that $|D \cap D'| = 1$, then compressing $\Sigma$ along $\pd N(\pd D \cup \pd D')$ yields a new surface $\Sigma'$, and it a simple exercise to show that $\Sigma'$ is a $(g-1,b)$-bridge surface for $K$ in $M$.  In this case, we say that $\Sigma$ is \emph{stabilized}. \\

If there are bridge disks $\Delta$ and $\Delta'$ for $(V,\A)$ and $(W,\n)$ such that $\Delta \cap \Delta'$ is a single point contained in $K$, there is an isotopy of $K$ along $\Delta \cup \Delta'$ resulting in a new bridge splitting $(M,K) = (V,\A') \cup_{\Sigma'} (W,\n')$, where each of $\A'$ and $\n'$ contains one fewer arc than $\A$ and $\n$, so that $\Sigma'$ is a $(g,b-1)$-bridge surface.  Here we say that $\Sigma$ is \emph{perturbed}. \\

Similarly, suppose $\Sigma$ is a $(g,1)$-surface, with bridge disks $\Delta$ and $\Delta'$ such that $\Delta \cap \Delta'$ is two points contained in $K$.  It follows that there is an isotopy which pushes $K$ onto a curve $c \subset \Sigma$, and we say $\Sigma$ is \emph{cancelable}.  If in addition there is a compressing disk $D$ for $\Sigma$ such that $D \cap c$ is a single point, we say that $\Sigma$ is \emph{removable}.  If $D \subset V$, we can push $c$ slightly into $V$, making $\Sigma$ a Heegaard surface for $E(K)$.  This procedure is carried out from the standpoint of fixing $\Sigma$ and isotoping $K$; equivalently, letting $D_{\Delta'} = \pd N(\Delta') \cap W$, compressing $\Sigma_K$ along the disk $D_{\Delta'}$ yields a closed surface $\Sigma'$, a Heegaard surface for $E(K)$. \\

There is one final case which allows us to simplify $\Sigma$:  If $D$ is a compressing disk for $(V,\A)$ and $C$ is a cut disk for $(W,\n)$ such that $|C \cap D|=1$, then compressing $\Sigma$ along $\pd N(\pd D \cup \pd C)$ yields a new surface $\Sigma'$, and it a simple exercise to show that $\Sigma'$ is a $(g-1,b+1)$-bridge surface for $K$ in $S^3$.  In this setting, we say that $\Sigma$ is \emph{meridionally stabilized}. \\

If $\Sigma$ is stabilized, perturbed, or meridionally stabilized, we say that $\Sigma$ is \emph{reducible}.  Otherwise, $\Sigma$ is \emph{irreducible}.  As with Heegaard surfaces, if there are compressing disks $D$ and $D'$ for $(V,\A)$ and $(W,\n)$ such that $D \cap D' = \emp$, we say that $\Sigma$ is \emph{weakly reducible}.  On the other hand, if every pair of compressing or bridge disks $D$ and $D'$ for $(V,\A)$ and $(W,\n)$ intersects away from the knot $K$, we say that $\Sigma$ is \emph{strongly irreducible}.  With these definitions, it is possible for an irreducible bridge surface $\Sigma$ to be neither weakly reducible nor strongly irreducible; we postpone the discussion of this case until Section \ref{esscan}. \\

In the literature there are variations of the above definitions using combinations of bridge disks, compressing disks, and cut disks.  The notions of weak reducibility and strong irreducibility above are the most natural in the context of this paper; however, we warn the reader that \emph{reducible} does not imply \emph{weakly reducible}.  More specifically, straightforward constructions show that a stabilized or perturbed bridge surface is weakly reducible, but a meridionally stabilized surface is not necessarily weakly reducible.  As an aside, we note that a meridionally stabilized surface is always $c$-weakly reducible, although to avoid confusion, we will not mention $c$-weak reducibility further.  For a detailed discussion of $c$-weak reducibility, see \cite{tomova}. \\



A properly embedded surface $S \subset E(K)$ is \emph{essential} if it is incompressible and not parallel into $\pd E(K)$.  Observe that since $\pd E(K)$ is a torus, if $S$ is essential, it must also be $\pd$-incompressible; that is, there is no embedded disk $\Delta$ such that $\pd \Delta$ is the endpoint union of arcs $\delta = \Delta \cap S$ and $\eps = \Delta \cap \pd E(K)$ with $\delta$ essential in $S$ (such a disk is called a \emph{$\pd$-compressing disk}).  The surface $S$ is \emph{meridional} if $S$ has nonempty boundary and $\pd S$ consists of meridian curves on $\pd E(K)$.

\section{Strongly irreducible bridge surfaces for torus knots}\label{strongly}

In this section, we show that a torus knot has a unique strongly irreducible bridge surface.  First, we require several more definitions, taken from \cite{zupan2}.  Let $M$ be an orientable, irreducible 3-manifold with boundary, $P \subset \pd M$ a subsurface, and $S \subset M$ a properly embedded surface.  A $P$-$\pd$-\emph{compressing} disk for $S$ is a $\pd$-compressing disk $\Delta$ for $S$ such that $\Delta \cap \pd M \subset P$.  We say that $S$ is $P$\emph{-essential} if $S$ is incompressible, there does not exist a $P$-$\pd$-compressing disk for $S$ in $M$, and $S$ is not a disk isotopic to a disk in $\pd M \setminus P$.  On the other hand, if $S$ is separating and admits compressing or $P$-$\pd$-compressing disks on either side but admits no pair of disjoint disks on opposite sides, we say that $S$ is $P$\emph{-strongly irreducible}.  As an example, a strongly irreducible bridge surface for a knot $K$ in $M$ is $\pd N(K)$-strongly irreducible in $E(K)$. \\

Now, we state a special case of Lemma 5.2 of \cite{zupan2}:

\begin{lemma}\label{L1}
Let $K$ be a knot in $S^3$, with $Q = \pd N(K)$ in $E(K)$.   Suppose $\Sigma$ is a strongly irreducible bridge surface for $(S^3,K)$, and let $A \subset E(K)$ be a properly embedded essential annulus.  Then either
\be
\item After isotopy, $\Sigma_K$ is transverse to $A$ and each component of $\Sigma_K \setminus \eta(A)$ is $Q$-essential in $E(K) \setminus \eta(A)$, or
\item After isotopy, $\Sigma_K$ is transverse to $A$, one component of $\Sigma_K \setminus \eta(A)$ is $Q$-strongly irreducible and all other components are $Q$-essential in $E(K) \setminus \eta(A)$.
\ee
\end{lemma}

Note that for a torus knot $K$, its exterior $E(K)$ contains precisely two essential surfaces: a Seifert surface $R$ and a properly embedded annulus $A$ such that $E(K) \setminus \eta(A) = V_1 \cup_A V_2$, where $V_1$ and $V_2$ are solid tori.  In addition, we may view each $V_i$ as the exterior of the unknot $K_0$.

\begin{lemma}\label{L2}
Suppose $\Sigma$ is a strongly irreducible bridge surface for a torus knot $K$, with $E(K) = V_1 \cup_A V_2$ as above.  Then for some $i$, $\Sigma_K \cap V_i$ is a bridge surface for the unknot $K_0$.
\begin{proof}
By Lemma \ref{L1}, for some component of $E(K) \setminus \eta(A)$, say $V_1$, we have that $\Sigma_K \cap V_1$ is $Q$-essential.  Since $\Sigma \cap K \neq \emp$, some component of $\Sigma_K \cap V_1$ cannot be a $\pd$-parallel annulus or disk, as these surfaces are either disjoint from $K$ or are $Q$-$\pd$-compressible.  Thus, some component of $\Sigma_K \cap V_1$ is a meridian disk intersecting $K$ either $p$ or $q$ times (and it follows that each component of $\Sigma_K \cap V_1$ must be this type of disk).  Viewing $V_1$ as a regular neighborhood of its unknotted core $K_0$, we see that $\Sigma_K \cap V_2$ is a bridge surface for $K_0$. \\
\end{proof}
\end{lemma}

Next, we invoke several well-known results about bridge surfaces for the unknot, proved by Otal \cite{otal} and Hayashi and Shimokawa \cite{hayshi}.

\begin{lemma}\label{L3}\cite{otal}
If $\Sigma$ is a $(0,b)$-bridge surface for the unknot with $b \geq 2$, then $\Sigma$ is perturbed. 
\end{lemma}

\begin{lemma}\label{L4}\cite{hayshi}
If $\Sigma$ is a $(g,b)$-bridge surface for the unknot with $g \geq 1$, then $\Sigma$ is perturbed or stabilized.
\end{lemma}

\begin{figure}[h!]
  \centering
    \includegraphics[width=1.0\textwidth]{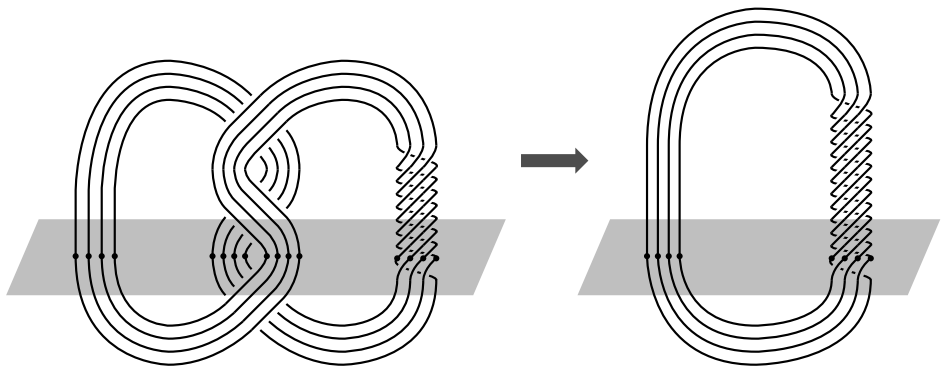}
    \caption{A $(13,4)$-cable of a $(0,2)$-splitting of the unknot (left) is the result of four perturbations of a cable of a $(0,1)$-splitting of the unknot (right).  The figure on the right is the unique strongly irreducible bridge surface for the $(13,4)$-torus knot.}
    \label{sisurf}
\end{figure}

Putting the three previous lemmas together, we arrive at the following conclusion:

\begin{lemma}\label{L5}
Suppose $\Sigma$ is a strongly irreducible $(g,b)$-bridge surface for a $(p,q)$-torus knot $K$.  Then $g=0$ and $b = \min\{p,q\}$.
\begin{proof}
Let $(S^3,K) = (V,\A) \cup_{\Sigma} (W,\n)$.  By Lemma \ref{L2}, we may suppose that $\Sigma_K \cap V_2$ is a $(g,b')$ bridge surface for $K_0$.  By Lemma \ref{L4}, if $g \geq 1$, then $\Sigma_K \cap V_2$ is perturbed or stabilized.  Since any perturbed or stabilized surface is weakly reducible, $\Sigma_K \cap V_2$ is weakly reducible; hence $\Sigma$ is also weakly reducible, a contradiction.  We remark that it is also a simple construction to show further that if $\Sigma_K \cap V_2$ is perturbed or stabilized, then $\Sigma$ is perturbed or stabilized (see Figure \ref{sisurf}). \\

If $g = 0$ and $b' \geq 2$, then Lemma \ref{L3} states that $\Sigma_K \cap V_2$ is perturbed, and thus $\Sigma$ is weakly reducible.  Finally, if $g=0$ and $b'=1$, then $\Sigma_K \cap V_2$ is the standard $(0,1)$-bridge surface for the unknot, and by Lemma \ref{L2}, $\Sigma$ is either a $(0,p)$- or $(0,q)$-surface.  However, if $\Sigma$ is a $(0,\max\{p,q\})$-surface, $\Sigma$ is easily seen to be perturbed (see, for example, \cite{ozawa}); hence $b = \min\{p,q\}$, completing the proof. \\
\end{proof}
\end{lemma}

It is a classical theorem of Schubert \cite{schub} (with a modern proof given by Schultens \cite{sch2}) that the bridge number of a $(p,q)$-torus knot is $\min\{p,q\}$.  In addition, by \cite{ozawa}, a $(p,q)$-torus knot $K$ has a unique $(0,\min\{p,q\})$-bridge surface (depicted at right in Figure \ref{sisurf}), which implies

\begin{corollary}\label{strong}
A torus knot $K$ has a unique strongly irreducible bridge surface; namely, a $\min\{p,q\}$-bridge sphere.
\end{corollary}

\section{Weakly reducible bridge surfaces for torus knots}\label{weak}

In this section, we turn our attention from strongly irreducible to weakly reducible bridge surfaces, showing that every weakly reducible bridge surface for a torus knot is reducible.  The following theorem of Tomova will be crucial:

\begin{theorem}\cite{tomova}\label{reducto}
Suppose that $\Sigma$ is a weakly reducible bridge surface for $(M,K)$.  Then one of the following occurs:
\bi
\item $\Sigma$ is reducible,
\item $\Sigma$ is removable, or
\item $E(K)$ contains an essential meridional or closed surface $S$.
\ei
\end{theorem}

For a torus knot $K \subset S^3$, its exterior $E(K)$ contains no essential meridional or closed surfaces; thus, if $\Sigma$ is a weakly reducible bridge surface for $(S^3,K)$, then $\Sigma$ is reducible or removable.  We devote the remainder of this section to proving that a removable bridge surface $\Sigma$ for a torus knot is reducible as well. \\

Suppose that $E(K) = H \cup_{\Sigma'} C$ is a Heegaard splitting, with $H$ a handlebody and $C$ a compression body.  We say that $\Sigma'$ is $\mu$-\emph{primitive} if there exists a compressing disk $D$ for $H$ and a vertical annulus $B$ in $C$ such that $\pd B = \gamma \cup m$, where $\gamma \subset \pd_+ C$ satisfies $|\gamma \cap \pd D| = 1$ and $m \subset \pd N(K)$ bounds a meridian of $N(K)$.  The next lemma uses this definition to help us understand removable bridge surfaces.

\begin{lemma}\label{L6}
Suppose that $(M,K) = (V,\A) \cup_{\Sigma} (W,\n)$ is a $(g,1)$-bridge splitting such that $\Sigma$ is removable, and let $\Sigma'$ be the corresponding Heegaard surface for $E(K)$.  If $\Sigma'$ is stabilized or $\mu$-primitive, then $\Sigma$ is stabilized.
\begin{proof}
We may view $\Sigma'$ as a Heegaard surface for $M$ bounding handlebodies $H_1$ and $H_2$, such that $K \subset H_1$ and $H_1 \setminus \eta(K)$ is a compression body.  In addition, $H_1$ contains an annulus $\Lambda$ with the property that $\pd \Lambda = K \cup \gamma$, where $\gamma \subset \Sigma'$.  Clearly, $\gamma$ is isotopic to $K$, and we choose $\Lambda$ so that $\gamma$ is the result of canceling the $(g,1)$-splitting given by $\Sigma$. \\

For any essential arc $\delta \subset \Lambda$, we may construct an isotopy which pushes $K$ across $\delta$ and into $H_2$ in a small neighborhood of $\delta \cap \gamma$.  This isotopy transforms $\Sigma'$ into a $(g,1)$-surface $\Sigma''$ which is equivalent to our original surface $\Sigma$. \\

Now, suppose $\Sigma'$ is stabilized, with compressing disks $D_1 \subset H_1 \setminus \eta(K)$ and $D_2 \subset H_2$ so that $|D_1 \cap D_2| = 1$.  Since $\Lambda$ is incompressible in $H_1 \setminus \eta(K)$, $D_1 \cap \Lambda$ cannot contain an essential simple closed curve, and there exists an essential arc $\delta \subset \Lambda$ with $\delta \cap D_1 = \emp$.  Using the construction above, we have that $D_1$ and $D_2$ avoid $K$; hence, $D_1$ and $D_2$ are compressing disks for $\Sigma''$ and $\Sigma''$ is stabilized. \\

On the other hand, if $\Sigma'$ is $\mu$-primitive, there is a vertical annulus $B \subset H_1 \setminus \eta(K)$ and compressing disk $D_2 \subset H_2$ such that $|B \cap D_2| = 1$.  Since $\pd B \cap \pd N(K)$ is meridional, $|\Lambda \cap B \cap \pd N(K)| = 1$, so $\Lambda \cap B$ contains a single essential arc $\delta$.  As above, there is an isotopy pushing $K$ across $\delta$ creating a $(g,1)$-surface $\Sigma''$.  This isotopy turns $B$ into a compressing disk for $\Sigma''$, and again $\Sigma''$ is stabilized.

\end{proof}
\end{lemma}

We may now invoke several known results:  

\begin{theorem}\cite{morsch}\label{TH1}
Let $K$ be a torus knot, and suppose $\Sigma'$ is a genus $g \geq 3$ Heegaard surface for $E(K)$.  Then $\Sigma'$ is stabilized.
\end{theorem}

\begin{theorem}\cite{brz, moriah}\label{TH2}
For a torus knot $K$, the exterior $E(K)$ has at most three distinct genus two Heegaard surfaces, up to isotopy.
\end{theorem}
We call these surfaces $\Sigma_u$, $\Sigma_l$, and $\Sigma_m$ (corresponding to upper, lower, and middle unknotting tunnels for $K$, respectively).  In certain cases, $\Sigma_m$ may be isotopic to either $\Sigma_l$ or $\Sigma_u$.  In any case, it is shown in \cite{moriah} that $\Sigma_u$ and $\Sigma_l$ are $\mu$-primitive. \\

Thus, if $\Sigma$ is a removable $(g,1)$-surface for a torus knot $K$, where $g \geq 3$ or the corresponding Heegaard surface $\Sigma'$ is isotopic to $\Sigma_u$ or $\Sigma_l$, then $\Sigma$ is stabilized by Lemma \ref{L6}.  To complete our classification of removable bridge surfaces, we need only consider the case in which $g = 2$ and $\Sigma' = \Sigma_m$.  The construction of $\Sigma_m$ is as follows: \\

\begin{figure}[h!]
  \centering
    \includegraphics[width=.6\textwidth]{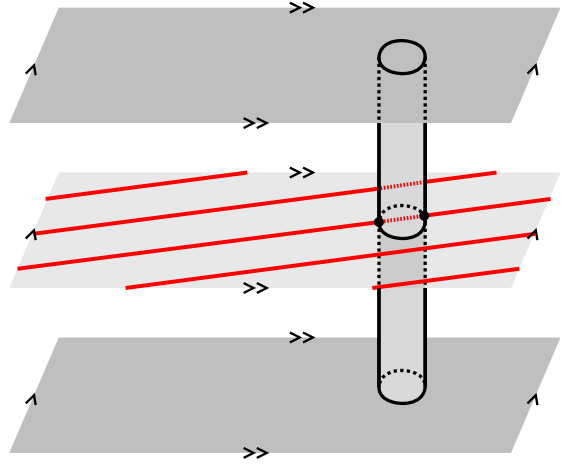}
    \caption{The $(2,1)$-bridge surface $\Sigma_*$ constructed from the Heegaard surface $\Sigma_m$.}
    \label{sigmam}
\end{figure}

Consider a Heegaard torus $T$ for $S^3$ with a collar neighborhood $T \X I$, and such that the $(p,q)$-torus knot $K$ is contained in $T \X \{1/2\}$.  Let $x \in T$ so that $\{x\} \X \{1/2\}$ misses $K$, and for $i=0,1$ let $V_i$ denote the solid torus disjoint from $T \X (0,1)$ bounded by $T \X \{i\}$.  Then $H = V_0 \cup V_1 \cup N(\{x\}) \X I$ is a genus 2 handlebody, and $C' = (T \setminus \eta(\{x\})) \X I$ is a genus 2 handlebody containing $K$, with $C=C' \setminus \eta(K)$ a compression body.  The surface $\Sigma_m$ is defined as $\pd H = \pd_+ C$.  There is an isotopy of $K$ supported in $\eta(T \X \{1/2\})$ which pushes an arc of $K$ into the tube $N(\{x\}) \X I$ to create a $(2,1)$-surface we call $\Sigma_*$.  See Figure \ref{sigmam}. \\

There is a bit of subtlety here:  Starting with the closed surface $\Sigma_m$, we may construct $\Sigma_*$ by pushing $K$ onto a curve $c_* \subset \Sigma_m$, followed by a perturbation of $c_*$.  Now, given any removable bridge surface $\Sigma$ which induces $\Sigma_m$, we may construct $\Sigma$ by pushing $K$ onto a curve $c \subset \Sigma_m$, followed by a perturbation of $c$.  However, a priori the curves $c$ and $c_*$ need not be isotopic in $\Sigma_m$, and we address this possibility in the next lemma.

\begin{lemma}\label{LL1}
If $\Sigma$ is a removable bridge surface which induces $\Sigma_m$, then $\Sigma$ is equivalent to $\Sigma_*$.
\begin{proof}
Let $E(K) = H \cup_{\Sigma_m} C$, with a curve $c_* \subset \Sigma_m$ whose perturbation yields $\Sigma_*$, as above, and let $c$ be a curve whose perturbation yields $\Sigma$.  Note that $c$ (resp. $c_*)$ cobounds an annulus $\Lambda$ (resp. $\Lambda_*$) with $K$ in the genus 2 handlebody $C' = C \cup N(K)$.  As in the proof Lemma \ref{L6}, an isotopy of $K$ along any essential arc $\delta \subset \Lambda$ (resp. $\delta_* \subset \Lambda_*$) yields a surface isotopic to $\Sigma$ (resp. $\Sigma_*$).  Thus, it suffices to show that there exist such essential arcs $\delta$ and $\delta_*$ which are isotopic. \\

There is a meridional annulus $A' \subset C$ that hits $c$ and $c_*$ exactly once; hence, $A'$ hits each of $\Lambda$ and $\Lambda_*$ in precisely one essential arc, say $\delta$ and $\delta^*$.  But $\delta$ and $\delta_*$ are also essential arcs in $A'$, and as such they are isotopic. \\
\end{proof}
\end{lemma}

\begin{lemma}\label{LL2}
The $(2,1)$-surface $\Sigma_*$ is stabilized.
\begin{proof}
Using the notation above, let $(S^3,K) = (V,\A) \cup_{\Sigma_*} (W,\n)$, where $V$ is the union of the Heegaard tori $V_0$ and $V_1$ vertically tubed across $T \X I$, $\A$ is an arc contained in $N(\{x\}) \X \{1/2\}$, $W = (T \setminus \eta(\{x\})) \X I$, and $\n$ is an essential arc in $(T \setminus \eta(\{x\})) \X \{1/2\}$.  Let $\Delta_W = \n \X [0,1/2]$, noting that $\Delta_W$ is a bridge disk for $(W,\n)$. \\

Let $D_V$ be a compressing disk for $(V,\A)$ contained in the solid torus $V_1$, and let $\gamma \subset (T \setminus \eta(\{x\})) \X \{1\}$ be an essential arc intersecting $\pd D_V$ exactly once.  Note that $D_W = \gamma \X I$ is a compressing disk in the handlebody $W$ (ignoring $\n)$, and let $D_* = \pd N(\Delta_W) \cap W$, so that $D_*$ is a compressing disk for $(W,\n)$, and $D_* \cap \Delta_W = \emp$. \\

If $D_W \cap D_* = \emp$, then $D_W$ avoids $\n$ and $\Sigma_*$ is stabilized.  If $D_W \cap D_* \neq \emp$, let $\delta$ be an arc of $D_* \cap D_W$ which is outermost in $D_*$.  A standard cut-and-paste argument provides a disk $D_W'$ such that $D_W' \cap (T \X \{1\}) = \gamma$, and $D_W'$ has fewer intersections with $D_*$ than $D_W$.  Finitely many repetitions yield such a disk $D_W'$ disjoint from $D_*$, completing the proof. \\
\end{proof}
\end{lemma}

We summarize this analysis with the following lemma:

\begin{lemma}\label{LL3}
If $\Sigma$ is a weakly reducible bridge surface for a torus knot $K$, then $\Sigma$ is reducible.
\begin{proof}
If $\Sigma$ is a weakly reducible $(g,b)$-bridge surface, then either it is reducible or removable by Theorem \ref{reducto}.  If $\Sigma$ is removable and $g \geq 3$, then $\Sigma$ is stabilized by Theorem \ref{TH1} and Lemma \ref{L6}.  Otherwise, $g = 2$, and the induced Heegaard splitting $\Sigma'$ is isotopic to either $\Sigma_u$, $\Sigma_l$, or $\Sigma_m$ by Theorem \ref{TH2}.  In the first two cases, $\Sigma$ is stabilized by Lemma \ref{L6}.  In the third case, Lemma \ref{LL1} yields that $\Sigma$ is equivalent to $\Sigma_*$, and $\Sigma_*$ is stabilized by Lemma \ref{LL2}, completing the proof. \\
\end{proof}
\end{lemma}

\section{Essentially cancelable bridge surfaces}\label{esscan}

Using our definitions, it is possible for an irreducible bridge surface $\Sigma$ to be neither weakly reducible nor strongly irreducible.  In this case, $\Sigma$ must be cancelable, and boundary compressing $\Sigma_K$ along a pair of canceling bridge disks yields an essential surface $S \subset E(K)$.  In this construction, $\pd S$ necessarily has longitudinal slope as a curve on $\pd N(K)$.  We call such a bridge surface \emph{essentially cancelable}. \\

\begin{figure}[h!]
  \centering
    \includegraphics[width=0.55\textwidth]{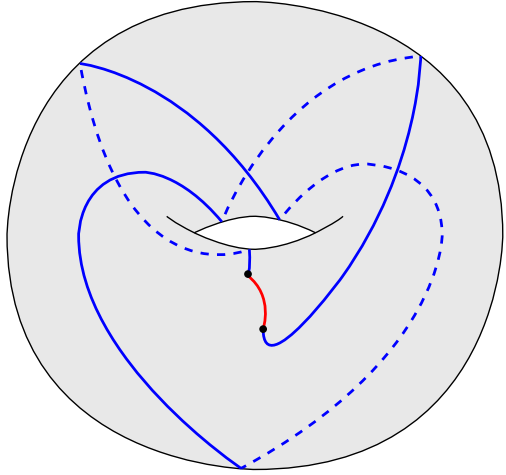}
    \caption{The $(1,1)$-surface $\Sigma_A$, where the blue and red arcs are on the inside and outside (resp.) of the Heegaard torus.}
    \label{tref}
\end{figure}

In this section, we consider essentially cancelable bridge surfaces for a $(p,q)$-torus knot $K$.  As noted above, the exterior $E(K)$ contains only two (connected) essential surfaces, the essential annulus $A$ discussed in Section \ref{strongly}, and a genus $(p-1)(q-1)/2$ Seifert surface $R$.  It is clear that the surface $S$ given in the construction above must be separating, and so there are two possibilities to consider:  Let $\Sigma_A$ denote the cancelable surface inducing $A$ (shown in Figure \ref{tref}), and let $\Sigma_R$ denote a cancelable surface inducing two copies of $R$.  We examine these cases separately. \\

It is clear that $\Sigma_A$ is a $(1,1)$-surface, and we have the following theorem of Morimoto:

\begin{theorem}\cite{morim}
A $(p,q)$-torus knot $K$ has a unique $(1,1)$-bridge surface.
\end{theorem}

As shown in \cite{schartom}, if $p=2$ or $q=2$, then $K$ is a 2-bridge knot and $\Sigma_A$ is meridionally stabilized.  Otherwise, the bridge number of $K$ is at least three, and $\Sigma_A$ is irreducible. \\

The last surface to classify is $\Sigma_R$, which we show falls into the weakly reducible case already examined in Section \ref{weak}.

\begin{lemma}\label{escan}
For a $(p,q)$-torus knot $K$, the surface $\Sigma_R$ is weakly reducible.
\begin{proof}
Let $\Sigma'$ denote a genus $(p-1)(q-1)$ Heegaard surface for $S^3$ such that $S^3 = V \cup_{\Sigma'} W$, $K \subset \Sigma'$, $\Sigma' \setminus \eta(K)$ is two copies of $R$.  It follows that any perturbation of $K$ will yield a $((p-1)(q-1),1)$-bridge surface equivalent to $\Sigma_R$. \\

Since $K$ is fibered, we may view $V$ as $R \X [0,1/2]$ and $W$ as $R \X [1/2,1]$.  Let $\varphi$ denote the monodromy of the fibration; that is, $E(K)$ is the mapping torus $(R \X I)/((x,0) \sim (\varphi(x),1))$.  It suffices to show that there exist disjoint essential arcs $\gamma_1,\gamma_2 \subset R$ such that $\gamma_1 \cap \varphi(\gamma_2) =\emp$ with the endpoints of $\gamma_2$ contained in the same component of $K \setminus \pd \gamma_1$.  In this case, we let $D_V = \gamma_1 \X [0,1/2]$ and $D_W = \gamma_2 \X [1/2,1]$, so that $D_V$ and $D_W$ are disjoint compressing disks for $V$ and $W$ each intersecting $K$ in two points.  After perturbation which pushes $K$ off of the two points of intersection with $D_V$ and the two points of intersection of $D_W$, the disks $D_V$ and $D_W$ become disjoint compressing disks for $\Sigma_R$, which is weakly reducible.  See Figure \ref{poff}.

\begin{figure}[h!]
  \centering
    \includegraphics[width=1.0\textwidth]{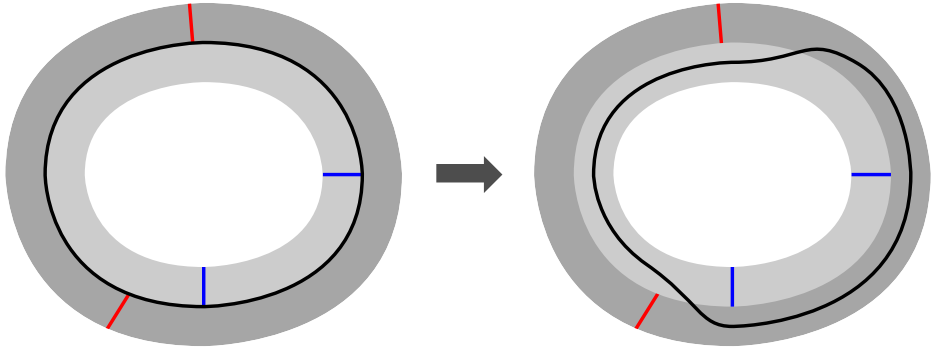}
    \caption{A small annular neighborhood of $K \subset \Sigma'$ transverse to $\Sigma'$ (left) and a perturbation of $K$ inducing $\Sigma_R$ (right).  Intersections with $D_V$ and $D_W$ are shown in red and blue, respectively.}
    \label{poff}
\end{figure}

To show the existence of such arcs $\gamma_1$ and $\gamma_2$, we use the following description of $R$:  Consider a Heegaard torus $T$ for $S^3$ containing $K$ and bounding solid tori $W_1$ and $W_2$.  The Seifert surface $R$ may be constructed by attaching $pq$ twisted bands from $p$ meridian disks of $W_1$ to $q$ meridian disks of $W_2$.  Suppose $p > q$ and let $D$ be a meridian disk of $W_2$ intersecting $K$ exactly $p$ times.  Choosing $x_1 \in D \cap K$ and an orientation for $K$ yields a natural ordering of the intersection points $x_1,x_2,\dots,x_p$ of $D \cap K$. \\

Now, let $\gamma_1$ be an arc in $D$ connecting $x_1$ and $x_2$, and let $\gamma_2$ be an arc connecting two other points $x_i$ and $x_j$ in $D \setminus \gamma_1$.  Observe that $\varphi(\gamma_2)$ is contained in another meridian disk component $D'$ of $R$; hence it is disjoint from $\gamma_1$.  We conclude that $\gamma_1$ and $\gamma_2$ satisfy the above conditions, completing the proof.

\end{proof}
\end{lemma}

\section{Conclusion}\label{conje}

In summary, we have
\begin{mtheorem}
Let $K$ be a torus knot which is not 2-bridge.  Then $K$ has a unique irreducible bridge surface of positive genus, a $(1,1)$-surface.
\begin{proof}
Suppose $\Sigma$ is an irreducible bridge surface for $K$.  Then $\Sigma$ is either strongly irreducible, weakly reducible, or essentially cancelable.  In the first case, $\Sigma$ is a minimal bridge sphere by Corollary \ref{strong}.  On the other hand, $\Sigma$ cannot be weakly reducible since $\Sigma$ would also be reducible by Lemma \ref{LL3}.  Finally, if $\Sigma$ is essentially cancelable, then $\Sigma$ is isotopic to either $\Sigma_A$ or $\Sigma_R$.  However, $\Sigma_R$ is reducible by Lemmas \ref{escan} and \ref{LL3}; thus $\Sigma_A$ is the unique irreducible bridge surface of positive genus. \\
\end{proof}
\end{mtheorem}

We conclude by remarking that in \cite{zupan}, the author shows that an $n$-fold iterated torus knot $K_n$ has a unique irreducible bridge sphere.  However, the author also shows in \cite{zupan2} (with an additional assumption on the cabling parameters) that for each genus $g=0,\dots,n+1$, the knot $K_n$ has an irreducible genus $g$ bridge surface.  This suggests the following conjecture: \\

\begin{conjecture}
Let $K_n$ be an $n$-fold iterated torus knot.  Then for each $g=0,\dots,n+1$, the knot $K_n$ has a unique irreducible genus $g$ bridge surface, and every $(g',b)$-surface with $g' > n+1$ is reducible.
\end{conjecture}

This paper settles the conjecture in the case that $n=0$.

\end{document}